\documentclass[11pt]{article}
\usepackage{cite}
\usepackage{mathrsfs}  
\usepackage{amssymb,amsfonts,amsmath}
\usepackage{enumerate}
\usepackage{indentfirst}
\usepackage{latexsym,bm}
\usepackage{algorithm}
\usepackage{makeidx}
\usepackage{fancybox}
\usepackage{multicol}  
\usepackage[latin1]{inputenc}
\usepackage{CJK}
\usepackage{color}     
\usepackage{pstricks,pst-node,pst-text,pst-3d}
\usepackage{tikz}
\usepackage{caption}
\usepackage{geometry}  
\makeatletter

\newcommand{\Rmnum}[1]{\expandafter\@slowromancap\romannumeral #1@}
\makeatother

\newtheorem{thm}{Theorem}[section]

\newtheorem{lem}[thm]{Lemma}

\newenvironment{pf}{{\noindent \it \bf Proof:}}{{\hfill$\Box$}\\}

\begin{document}
\begin{CJK}{GBK}{song}
\title{\bf The normalized Laplacian and related indexes of graphs with edges blew up by cliques}
\author{
\small Qi Ma, Zemin Jin\footnote{Corresponding author.  Email: {\tt silvester.ma@outlook.com (Ma), zeminjin@zjnu.cn (Jin)}} \\
\small   Department of Mathematics, Zhejiang Normal
University\\
\small Jinhua 321004, P.R. China \\
}

\date{}
\maketitle

\begin{abstract}
In this paper, we introduce the clique-blew up graph $CL(G)$ of a given graph $G$, which is obtained from $G$ by replacing each edge of $G$ with a complete graph $K_n$. We characterize all the normalized Laplacian spectrum of the grpah $CL(G)$ in term of the given graph $G$. Based on the spectrum obtained, the formulae to calculate the multiplicative degree-Kirchhoff index, the Kemeny's constant and the number of spanning trees of $CL(G)$  are derived well. Finally, the spectrum and indexes of the clique-blew up iterative graphs are present.
\\
[2mm]{\bf Key Words}: adjacent matrix; normalized Laplacian; multiplicative degree-Kirchhoff index; Kemeny's constant; spanning tree. \\
[2mm] {\bf AMS subject classification (2010)}:  05C50, 05C76.
\end{abstract}


\section{Introduction}

\subsection{Notions and definitions}
We consider a simple and connected graph $G=(V(G),E(G))$ with $n$ vertices and denote the vertex set of $G$ by $V(G)=\{1, 2, \cdots, n\}$. For any two adjacent vertices $i$ and $j$, we denote it by $i\sim j$. Denote the degree of a vertex $i$  by $d_i$ in $G$.  Let $A_G$ be the {\it adjacency matrix} of $G$, where the $(i,j)$-entry equals to $1$ if  $i\sim j$  and
$0$ otherwise. Clearly $A_G$ is an $n\times n$ matrix. Let $D_G=diag(d_1, d_2, \cdots, d_n)$ is the {\it diagonal matrix} of vertex degrees of $G$, where $d_i$ is the degree of $i$ in $G$. The matrix $L_G=D_G-A_G$ is called the {\it Laplacian matrix} of $G$.  Given a matrix $M$,  let $M(i,j)$ denote the $(i,j)$-entry of $M$. For the eigenvalue $\lambda$ of the matrix $M$, denote by  $m_{M} (\lambda)$ the multiplicity of $\lambda$ in $M$.

Given a graph  $G$, one can always define the random walk on $G$ as a Markov chain $X_n, n\geq0$.  The probability of  jumping from the current vertex $i$  to another vertex $j$ is $p_{ij}$, where $p_{ij}=\frac{1}{d_i}$ if $i$ and $j$ are adjacent and $p_{ij}=0$ otherwise, i.e.,
\begin{equation*}
p_{ij}=
\begin{cases}
\frac{1}{d_i}, & if \ i\sim j,  \\
{0}, &otherwise.
\end{cases}
\end{equation*}
The matirx $P_G=(p_{ij})_{n\times n}$ is the {\it transition probability matrix} for the random walk defined on $G$. It is clear that $P_G=D_{G}^{-1}A_G$.  The {\it normalized Laplacian matrix} of the graph
$G$ is defined to be
\begin{equation*}
\mathcal{L}_G=I-D_{G}^{\frac{1}{2}}P_{G}D_{G}^{-\frac{1}{2}},
\end{equation*}
where $I$ is an $n\times n$ identity marix.  Let $\delta _{ij}$ be the {\it Kronecker delta}, where $\delta _{ij}=1$ if $i=j$ and $\delta _{ij}=0$ otherwise.  According to the definition of $\mathcal{L}_G$,  we have that:
\begin{equation*}
\mathcal{L}_G(i,  j)=\delta _{ij} - \frac{A_G(i,  j)}{\sqrt {d_i d_j}}.
\end{equation*}
  The eigenvalues of $\mathcal{L}_G$ are non-negative because $\mathcal{L}_G$ is Hermitian to
$I-P_G=D_{G}^{-1}L_G$. For the $n$ eigenvalues of $\mathcal{L}_G$, we label them by $\lambda_1\leq \lambda_2\leq\cdots\leq\lambda_n$. Define the  {\it normalized Laplacian spectrum} on  $\mathcal{L}_G$ of the graph  $G$  as
$\sigma=\{\lambda_1,  \lambda_2,  \ldots,  \lambda_n\}$.

Often the normalized Laplacian spectrum of graphs can be used to characterize parameters of graphs, see \cite{chung}. Recently, one of very interesting applications of the spectrum of graphs is to study the the electric network.  Klein and Randi\'{c}\cite{klein}  proposed  a new distance function called {\it resistance distance} between two vertices in graphs.  Assume that there is a  unit resistor on every edge of the graph $G$.  When we attach  a
battery at two vertices $i$  and $j$,  the  resistance distance between $i$ and $j$, denoted by $r_{ij}$, is the electrical resistance between $i$  and $j$ in $G$.  For more recent results about resistance distances, one
can refer to \cite{jkk,somodi}.  Chen and Zhang \cite{chenzhang} proposed a new index called  the {\it multiplicative degree-Kirchhoff index} (see\cite{fenggut}) which is  defined as $
Kf^{*}(G)=\sum\limits_{i<j}d_{i}d_{j}r_{ij}.
$
There is a close relationship between the multiplicative degree-Kirchhoff index  and  the  normalized spectrum. In recent years, more and more results relating to  the normalized Laplacian spectrum and the multiplicative degree-Kirchhoff index of some  graphs have been obtained, see \cite{chenjost,fenggut,fengyu,huangli,huangliq,huangzhou,
liu,xie1,xie2}.  The {\it Kemeny's constant} $K_e(G)$ of $G$ is defined as the expected number of steps  for the transition from an initial vertex $i$ to a target vertex $j$, which can be selected randomly according to a stationary distribution of unbiased random walks on $G$. The Kemeny's constant provides an interesting and novel quantity for finite  ergodic Markov chains, which is unrelated to the initial state of the  Markov chain
\cite{levene,hunter}.

\subsection{Backgrounds}

 Many graph invariants, including the multiplicative degree-Kirchhoff index,  the Kemeny's constant, the number of spanning trees, can be calculated in term of the spectrum of the graph. In recent years, some researchers focused on blowing up all the edges of a given graph by replacing each edge with some another graph. The spectrum of the resulting graph always can be characterized in term of the given graph.

  Xie et al. \cite{xie1} initially replaced each edge of a graph $G$ with a triangle. They added a parallel path of lengths two between each two adjacent vertices. The spectrum of the normalized Laplacian of the new graph are characterized in term of $G$.  Later,
  Wang et al.  \cite{wang}  generalized the result of  \cite{xie1} by replacing each edge with $k$ triangles, i.e., they added $k$ edge-disjoint paths of length two between each two adjacent two vertices.  Li and  Hou \cite{li} blew up each edge of $G$ to a $4$-cycle by adding a new path of length three between each two adjacent vertices. The resulting graph is called the quadrilateral graph $Q(G)$. Huang and  Li \cite{huanglik} further added $k$  paths of length three between each two adjacent vertices to get the so-called $k$-quadrilateral graph $Q^k(G)$ of $G$. Luckily, the normalized Laplacian spectra of these resulting graphs can be characterized completely in term of the initial graph $G$.  As applications, one can calculate the multiplicative degree-Kirchhoff index, the Kemeny's constant and the number of spanning trees of of these graphs again in term of the initial graph.

   Pan  et al. \cite{panyingui} introduced an analogue method to replace the edges of a graph. They added  a triangle or a $4$-cycle between each two adjacent vertices and connected these vertices in a suitable way. The subdivision graph was considered in \cite{xie2}. More ideas to blow up  the edges of a given graph were studied in \cite{das,huangli}.  The authors \cite{heli,lili,scli,maxiaolin} considered the graph chains, which are obtained by replacing only one edge of the given graph iteratively with some special structures.  In addition to the spectra of the obtained graphs above, the authors \cite{chen,huang,wang,zhu}  studied the hitting times of the random walks on these graphs.

\section{Preliminaries}

\begin{figure}[H]
\begin{center}
\begin{tikzpicture}

\fill (-0.5,0) circle(2pt);
\fill (0.5,0) circle(2pt);
\fill (0,0.866) circle(2pt);

\fill (2.5,0.866) circle(2pt);
\fill (3,0) circle(2pt);
\fill (2,0) circle(2pt);

\fill (3,1.732) circle(2pt);
\fill (3.5,0.866) circle(2pt);
\fill (3.5,0) circle(2pt);

\fill (2,1.732) circle(2pt);
\fill (1.5,0.866) circle(2pt);
\fill (1.5,0) circle(2pt);

\fill (2.5,-1) circle(2pt);
\fill (1.5,-0.866) circle(2pt);
\fill (3.5,-0.866) circle(2pt);

\draw(2.5,0.866)--(3,0);
\draw(2.5,0.866)--(2,0);
\draw(2,0)--(3,0);

\draw(3,1.732)--(2.5,0.866);
\draw(3.5,0.866)--(2.5,0.866);
\draw(3.5,0)--(2.5,0.866);

\draw(3,1.732)--(3,0);
\draw(3.5,0.866)--(3,0);
\draw(3.5,0)--(3,0);

\draw(3,1.732)--(3.5,0.866);
\draw(3,1.732)--(3.5,0);
\draw(3.5,0.866)--(3.5,0);

\draw(2,1.732)--(2.5,0.866);
\draw(1.5,0.866)--(2.5,0.866);
\draw(1.5,0)--(2.5,0.866);

\draw(2,1.732)--(2,0);
\draw(1.5,0.866)--(2,0);
\draw(1.5,0)--(2,0);

\draw(2,1.732)--(1.5,0.866);
\draw(2,1.732)--(1.5,0);
\draw(1.5,0.866)--(1.5,0);

\draw(2.5,-1)--(2,0);
\draw(1.5,-0.866)--(2,0);
\draw(3.5,-0.866)--(2,0);

\draw(2.5,-1)--(3,0);
\draw(1.5,-0.866)--(3,0);
\draw(3.5,-0.866)--(3,0);

\draw(2.5,-1)--(1.5,-0.866);
\draw(2.5,-1)--(3.5,-0.866);
\draw(3.5,-0.866)--(1.5,-0.866);

\draw(-0.5,0)--(0,0.866);
\draw(-0.5,0)--(0.5,0);
\draw(0.5,0)--(0,0.866);
\end{tikzpicture}
\end{center}
\caption{The graph $G=K_3$ and its clique-blew up graph $CL(G)$ for $n=5$.}\label{ck3}
\end{figure}
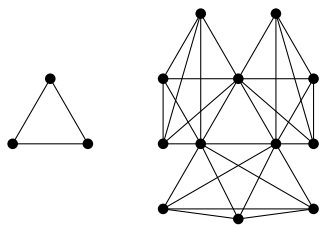

 Throughout all the paper, let $n\geq 3$ and $G$ be a simple and connected graph with $N_0$ vertices and $E_0$ edges.   For any edge $e$,   we add $n-2$ vertices, $k^e_i,i=1, 2, \cdots , {n-2}$, so
that all these vertices together with the end-vertices of $e$  form a $K_n$.  The resulting graph is called the {\it clique-blew up graph} and  written by $CL(G)$.  The Figure \ref{ck3} gives an example of the clique-blew up graph for
$G=K_3$ and $n=5$.

Let $V_{N}$ be the set of all the newly added vertices in $CL(G)$ and $V_{O}$ be the set of the vertices inherited from $G$.   That is,  the vertex set $V(CL(G))$ of $CL(G)$ is the union of $V_{N}$ and $V_{O}$.   We denote by $N_1$ the total number of vertices and $E_1$  the total number of edges of $CL(G)$.   It  is clear  that
 $E_1=\frac{n(n-1)E_0}{2}$   and  $N_1=N_0+(n-2)E_0.$

\begin{lem}\label{sp} {\upshape \cite{chung}}
Let $G$ be a connected graph with $n$ vertices ,  and $\mathcal{L}_G$ be the  normalized Laplacian matrix of  $G$.  The normalized Laplacian spectrum of $G$ is $\sigma=\{0=\lambda_1<\lambda_2\leq\cdots\leq\lambda_n\}$.  We
have

{\upshape (\romannumeral1)}  $\frac{n}{n-1}\leq\lambda_n\leq2$ with $\lambda_n=2$ if and only if $G$ is bipartite;

{\upshape (\romannumeral2)} If $G$ is bipartite, then for  any eigenvalue $\lambda_i$ of $\mathcal{L}_G$,   $2-\lambda_i$ is also an eigenvalue of $\mathcal{L}_G$ and $m_{\mathcal{L}_G}
(\lambda_i)=m_{\mathcal{L}_G}(2-\lambda_i)$.
\end{lem}

By determining the spectrum on the normalized Laplacian of $G$,   the  specific calculation  formulae  of the multiplicative degree-Kirchhoff  index,  the  Kmemeny's constant and the  number  of spanning trees of graph $G$
can  be  listed as follows.
\begin{lem}\label{lm1.2}
Let $G$ be a connected graph with $n$ vertices and $m$ edges,  $\sigma=\{0=\lambda_1<\lambda_2\leq\cdots\leq\lambda_n\}$ is the spectrum on the normalized Laplacian $\mathcal{L}_G$ of $G$.  Then

{\upshape (\romannumeral1)\cite{chenzhang}} The multiplicative degree-Kirchhoff  index of $G$ is
$
Kf^{*}(G)=2m\sum\limits_{i=2}^{n}\frac{1}{\lambda_i}.
$

{\upshape (\romannumeral2)\cite{butler}} The Kemeny's constant of $G$ is
$
K_e(G)=\sum\limits_{i=2}^{n}\frac{1}{\lambda_i}.
$

{\upshape (\romannumeral3)\cite{chung}}  The number $\tau (G)$ of spanning trees of $G$ is
$
\tau (G)=\frac{1}{2m}\prod\limits_{i=1}^{n}d_i\cdot \prod\limits_{k=2}^{n}\lambda_k.
$

{\upshape (\romannumeral4)} From {\upshape (i) and (ii)}, we have
$Kf^{*}(G)=2mK_e(G).
$
\end{lem}


\section{The normalized Laplacian spectrum of $CL(G)$}\label{nlsp}
 For the clique-blew up graph $CL(G)$ of $G$,   the normalized Laplacian of $CL(G)$ is denoted by $\mathcal{L}_{C}$. Denote the degree of the vertex $i\in{V(CL(G))}$ by $d_{i}^{'}$.  Let $A_{C}$ be the {\it adjacency matrix} of $CL(G)$ and $D_{C}$ be the degree matrix of $CL(G)$.  Let $N_{G}=D_G^{-\frac{1}{2}}A_{G}D_{G}^{-\frac{1}{2}}$ and $N_{C}=D_{C}^{-\frac{1}{2}}A_{C}D_{C}^{-\frac{1}{2}}$.  For the incidence matrix of a connected graph, we have the following result.


\begin{lem}\label{bx}\upshape{\cite{cvetkovic}}
Let $B$ be the  incidence matrix of a connected graph $G$ with $n$ vertices.  Then
\begin{equation*}
rank(B)=
\begin{cases}
{n-1}, &\text{if {$G$} is bipartite}, \\
{n}, &\text{if {$G$} is non-bipartite}.
\end{cases}
\end{equation*}
\end{lem}\vskip2mm

At first, we consider the eigenvalue and its eigenvector in the graph $CL(G)$.  Let $v=({v}_{1},  {v}_{2},  \cdots, {v}_{N_1})^{T}$ be an eigenvector with respect to the eigenvalue $\lambda$ of $\mathcal{L}_{C}$,  i.e.,
\begin{equation}\label{lqv}
\mathcal{L}_{C}v=(I-N_{C})v=\lambda v.
\end{equation}
For any vertex $u\in V(CL(G))$,  the Eq. (\ref{lqv}) indicates that
\begin{equation}\label{1-lamb}
(1-\lambda)v_{u}=\sum\limits_{p=1}^{N_1}N_{C}(u,  p)v_{p}=\sum\limits_{p=1}^{N_1}\frac{A_{C}(u,  p)}{\sqrt{d_{u}^{'}d_{p}^{'}} }v_p.
\end{equation}
For any vertex $i\in V_{O}$, denote by $N_{O}(i)$ the set of neighbors of $i$ in $G$. Let $e$ be an edge with end vertices $i$ and $j$ in $G$.  By the construction of
$CL(G)$ and Eq. (\ref{1-lamb}),   we have
\begin{equation}\label{1-lamvi'}
\begin{aligned}
(1-\lambda)v_{i}&=\sum\limits_{j\in N_{O}(i)}\frac{v_{j}}{\sqrt {d_{i}^{'}d_{j}^{'}}}+\sum\limits_{e\in E(G) \ \mbox{is incident with }\ i}\sum\limits_{l=1}^{n-2}\frac{v_{k^e_l}}{\sqrt {d_{i}^{'} d_{k^e_l}^{'}}}\\
&=\sum\limits_{j\in N_{O}(i)}\frac{v_{j}}{(n-1)\sqrt {d_{i}d_{j}}}+\sum\limits_{e\in E(G) \ \mbox{is incident with }\ i}\sum\limits_{l=1}^{n-2}\frac{v_{k^e_l}}{(n-1)\sqrt {d_{i}}}.
\end{aligned}
\end{equation}
Similarly,  for the new vertices $v_{k^e_1}$ and $v_{k^e_2}$ corresponding to the edge  $e\in E(G)$ with end vertices $i$ and $j$,  we have
\begin{small}
\begin{equation}\label{vk1}
\begin{aligned}
(1-\lambda)v_{k^e_1}&=\frac{v_{i}}{\sqrt{d_{i}^{'}d_{k^e_1}^{'}}}
+\frac{v_{j}}{\sqrt{d_{j}^{'}{d^{'}_{k^e_1}}}}+\frac{v_{k^e_2}}{\sqrt {d_{k^e_1}^{'} d_{k^e_2}^{'}}}+\frac{v_{k^e_3}}{\sqrt {d_{k^e_1}^{'} d_{k^e_3}^{'}}}+\cdots+\frac{v_{k^e_{n-2}}}{\sqrt {d_{k^e_1}^{'}
d^{'}_{k^e_{n-2}}}}\\
&=\frac{v_{i}}{(n-1)\sqrt{d_{i}}}
+\frac{v_{j}}{(n-1)\sqrt{d_{j}}}+\frac{v_{k^e_2}
+v_{k^e_3}+\cdots+v_{k^e_{n-2}}}{n-1},
\end{aligned}
\end{equation}
\end{small}
\begin{small}
\begin{equation}\label{vk2}
\begin{aligned}
(1-\lambda)v_{k^e_2}&=\frac{v_{i}}{\sqrt{d_{i}^{'}d_{k^e_2}^{'}}}
+\frac{v_{j}}{\sqrt{d_{j}^{'}{d^{'}_{k^e_2}}}}+\frac{v_{k^e_1}}{\sqrt {d_{k^e_2}^{'} d_{k^e_1}^{'}}}+\frac{v_{k^e_3}}{\sqrt {d_{k^e_2}^{'} d_{k^e_3}^{'}}}+\cdots+\frac{v_{k^e_{n-2}}}{\sqrt {d_{k^e_2}^{'}
d^{'}_{k^e_{n-2}}}}\\
&=\frac{v_{i}}{(n-1)\sqrt{d_{i}}}
+\frac{v_{j}}{(n-1)\sqrt{d_{j}}}+\frac{v_{k^e_1}
+v_{k^e_3}+\cdots+v_{k^e_{n-2}}}{n-1}.
\end{aligned}
\end{equation}
\end{small}

The following lemma shows the relationship between the normalized Laplacian eigenvalues of $CL(G)$ and $G$.
\begin{lem}\label{(n-1)times}
 Let $\lambda$ be an eigenvalue of $\mathcal{L}_{C}$ such that $\lambda \neq  {n\over n-1}$ and ${2\over n-1}$.  Then $(n-1)\lambda$ is an eigenvalue of $\mathcal{L}_{G}$ with $m_{\mathcal{L}_{C}}(\lambda)=m_{\mathcal{L}_{G}}((n-1)\lambda)$.
\end{lem}
\begin{pf}  Let $v=({v}_{1},  {v}_{2},  \cdots, {v}_{N_1})^{T}$ be an eigenvector with respect to the eigenvalue $\lambda$ of $\mathcal{L}_{C}$.  Let $e\in E(G)$ with end vertices $i$ and $j$.   Since $\lambda \neq  {n\over n-1}$, from  Eqs. \eqref{vk1} and \eqref{vk2}, we have $v_{k^e_1}=v_{k^e_2}$. For the same reason, we  can easily get
\begin{equation}\label{vk'}
v_{k^e_1}=v_{k^e_2}=\cdots=v_{k^e_{n-2}}.
\end{equation}
For convenience, let  $v_{k^e_1}=x_e$.
Substituting Eq. \eqref{vk'} into Eqs. \eqref{1-lamvi'} and \eqref{vk1}, we have
\begin{equation}\label{1-lamvi}
(1-\lambda)v_{i}=\sum\limits_{j\in N_{O}(i)}\frac{v_{j}}{(n-1)\sqrt {d_{i}d_{j}}}+\sum\limits_{e\in E(G) \ \mbox{is incident with }\ i}\frac{n-2}{(n-1)\sqrt {d_{i}}}x_{e},
\end{equation}
\begin{equation}\label{1-lamvk}
(2-(n-1)\lambda)x_{e}=\frac{v_{i}}{\sqrt{d_{i}}}+
\frac{v_{j}}{\sqrt{d_{j}}}.
\end{equation}
Combining Eqs.  (\ref{1-lamvi}) and (\ref{1-lamvk}),   for $\lambda\neq {{n\over {n-1}}\ \mbox{and} \  \frac{2}{n-1}}$,   it follows that
\begin{equation*}
(1-\lambda)v_i=\frac{n-2}{(n-1)(2-n\lambda+\lambda)}v_i+\sum\limits_{j\in N_{O}(i)}\frac{n-n\lambda+\lambda}{(n-1)(2-n\lambda+\lambda)\sqrt{d_i d_j}}v_j,
\end{equation*}
i.e.,
\begin{equation}\label{1-(n-1)lam}
(1-(n-1)\lambda)v_i=\sum\limits_{j\in N_{O}(i)}\frac{v_j}{\sqrt{d_i d_j}}.
\end{equation}
holds for $\lambda\neq   {n\over n-1}$ and ${2\over n-1}$.

From Eq.  (\ref{1-(n-1)lam}),   it is obvious that $1-(n-1)\lambda$ is an eigenvalue of the matrix $N_G$ for $\lambda\neq    {n\over n-1}$ and ${2\over n-1}$.   So for any eigenvalue  $\lambda$ ($\lambda\neq   {n\over
n-1}$ and ${2\over n-1}$) and a corresponding eigenvector $v$ of $\mathcal{L}_C$,    $(n-1)\lambda$ and  $(v_i)^T_{i\in V_{O}}$ are an eigenvalue  and a corresponding eigenvector of $\mathcal{L}_G$, respectively.  This
implies that $m_{\mathcal{L}_G}((n-1)\lambda)\geq  m_{\mathcal{L}_C}(\lambda)$.

 On the other hand, for any eigenvalue $(n-1)\lambda $ ($(n-1)\lambda\neq 0,2$) and a corresponding eigenvector  $(v_i)^T_{i\in V_{O}}$ of $\mathcal{L}_G$,  the value  $\lambda$ is an eigenvalue of $\mathcal{L}_C$. Also, the vector determined by
 $(v_i)^T_{i\in V_{O}}$  and Eq. (\ref{1-lamvk}) is a corresponding eigenvector for the eigenvalue $\lambda$ of $\mathcal{L}_C$. Hence  $m_{\mathcal{L}_G}((n-1)\lambda)\leq  m_{\mathcal{L}_C}(\lambda)$.
 So we have that $m_{\mathcal{L}_G}((n-1)\lambda)= m_{\mathcal{L}_C}(\lambda)$.
The proof is completed.
\end{pf}

Now we give a complete representation about the normalized Laplacian eigenvalues and corresponding eigenvectors of $CL(G)$ as follows.
\begin{thm}\label{sp-Q}
Let $G$ be a simple  connected graph with $N_0$  vertices and $E_0$ edges and $CL(G)$ be the clique-blew up  graph of $G$.  The normalized Laplacian spectrum of $CL(G)$ can  be obtained as following

{\upshape (\romannumeral1)} The value 0 is an eigenvalue of $\mathcal{L}_C$ with the multiplicity 1;

{\upshape (\romannumeral2)} If $\lambda\ (\lambda \neq 0,2)$ is an eigenvalue of $\mathcal{L}_G$,  then the value $\frac{\lambda}{n-1}$ is an eigenvalue of $\mathcal{L}_C$ and $m_{\mathcal{L}_C}(\frac{\lambda}{n-1})=m_{\mathcal{L}_G}(\lambda)$;

{\upshape (\romannumeral3)} If $G$ is non-bipartite, then the value $\frac{2}{n-1}$ is an eigenvalue of $\mathcal{L}_C$ with the multiplicity $E_0-N_0$;

{\upshape (\romannumeral4)} If $G$  is bipartite,  then the value $\frac{2}{n-1}$ is an eigenvalue of $\mathcal{L}_C$ with the multiplicity $E_0-N_0+1$;

{\upshape (\romannumeral5)} The value $\frac{n}{n-1}$ is the eigenvalue of $\mathcal{L}_C$ with the multiplicity
$(n-3)E_0+N_0$.
\end{thm}

\begin{pf}
{\upshape (\romannumeral1)}  It is obvious from Lemma \ref{sp}.

{\upshape (\romannumeral2)}  It follows from  Lemma \ref{(n-1)times} that the statement holds obviously.\vskip2mm

 Since each eigenvalue $\lambda$ ($\lambda\neq \frac{2}{n-1},\frac{n}{n-1}$) of $\mathcal{L}_C$ and its multiplicity  have been determined in the statements above, here we only need to consider the eigenvalue $\lambda\in \{\frac{2}{n-1}, \frac{n}{n-1}\}$.\vskip2mm

Let $v=({v}_{1},  {v}_{2},  \cdots, {v}_{N_1})^{T}$ be an eigenvector with respect to the eigenvalue $\lambda={2\over n-1}$ of $\mathcal{L}_{C}$. Let $i\in V_{O}$ and $e\in E(G)$ with end vertices $i$ and $j$. For $n\geq 3$, from  Eqs. \eqref{vk1} and \eqref{vk2}, we have $\frac{n-2}{n-1}(v_{k_e^1}-v_{k_e^2})=0$, that is to say, $v_{k^e_1}=v_{k^e_2}$. For the same reason, we  can easily get
\begin{equation}\label{vk'.}
v_{k^e_1}=v_{k^e_2}=\cdots=v_{k^e_{n-2}}.
\end{equation}
For convenience, let  $v_{k^e_1}=x_e$.
Substituting Eq. \eqref{vk'.} and $\lambda=\frac{2}{n-1}$ into  Eqs. \eqref{1-lamvi'} and \eqref{vk1}, we have that
\begin{equation}\label{1-lamvi.}
\frac{n-3}{n-1}v_{i}=\sum\limits_{j\in N_{O}(i)}\frac{v_{j}}{(n-1)\sqrt {d_{i}d_{j}}}+\sum\limits_{e\in E(G) \ \mbox{is incident with }\ i}\frac{n-2}{(n-1)\sqrt {d_{i}}}x_{e}
\end{equation}
and \begin{equation}\label{vivj}
\frac{v_i}{\sqrt{d_i}}=-\frac{v_j}{\sqrt{d_j}}.
\end{equation}
\vskip2mm

\upshape{(\romannumeral3)} Let $G$ be non-bipartite. Suppose that $C$ is an odd cycle in $G$ of length $l$ with its vertices $i_1,  i_2,  .  .  .  ,  i_l$ in turn.  By Eq.  (\ref{vivj}), we have  $$\frac{v_{i_1}}{\sqrt {d_{i_1}}}=-\frac{v_{i_2}}{\sqrt
{d_{i_2}}}=\frac{v_{i_3}}{\sqrt {d_{i_3}}}=\cdots=\frac{v_{i_l}}{\sqrt{ d_{i_l}}}=-\frac{v_{i_1}}{\sqrt{ d_{i_1}}},$$
      which implies that $v_{i_k}=0,   k=1,  2,  \cdots,  l$.   Since $G$ is connected,   it holds that
\begin{equation}\label{vi}
v_i=0\  \mbox{for all}\ i\in V_{O}.
\end{equation}
Together with Eq. \eqref{1-lamvi.},  we have that for $i\in V_{O}$
\begin{equation}\label{sumvi}
\sum\limits_{e\in E(G) \ \mbox{is incident with}\ i}x_e=0.
\end{equation}

Therefore,   the eigenvectors $v=(v_1,  v_2,  .  .  .  ,  v_{N_1})^{T}$ associated with $\lambda=\frac{2}{n-1}$ can be determined by Eqs.  (\ref{vk'.})(\ref{vi}) and (\ref{sumvi}).   Notice that
$v_{k^e_1}=v_{k^e_2}=\cdots=v_{k^e_{n-2}}=x_e$.  Let $\mbox{\bf x}=(x_e)^T$ which is an $E_0$ dimensional vector by the construction of $CL(G)$.  It is easy to see that Eq. (\ref{sumvi}) is equivalent to the equation system
$B\mbox{\bf x}=0$,   where $B$ is the incident matrix of $G$.   By Lemma \ref{bx},   for $B\mbox{\bf x}=0$,   the number of solutions in its basic solution system is $E_0-N_0$ when $G$ is non-bipartite,   i.e.,
$m_{\mathcal{L}_C}(\frac{2}{n-1})=E_0-N_0$.
\vskip2mm

\upshape{(\romannumeral4)} Let $G$ be bipartite. Substituting Eq. \eqref{vivj} into Eq. \eqref{1-lamvi.}, we have that for $ i\in V_{O}$,
\begin{equation}\label{vixe}
\sqrt{d_i}v_i= \sum\limits_{e\in E(G) \ \mbox{is incident with }\ i}x_e.
\end{equation}
Let
$
\frac{v_1}{\sqrt{d_1}}=t.
$
 Denote by  $X$ and $Y$ the partite sets of the graph $G$ and without loss of generality, let $1\in X$. Then from Eq. \eqref{vivj}, we have that  $\frac{v_i}{\sqrt{d_i}}=t$ if $i\in X$, and $\frac{v_i}{\sqrt{d_i}}=-t$ if
 $i\in Y$.
 According to Eq. \eqref{vixe}, we have that for each $i\in V_O$,
 \begin{equation}\label{xevit}
 \begin{aligned}
&\sum\limits_{e\in E(G) \ \mbox{is incident with }\ i}x_e-d_i t=0 \ \mbox{if}\ i\in X,\\
&\sum\limits_{e\in E(G) \ \mbox{is incident with }\ i}x_e+d_i t=0 \ \mbox{if}\ i\in Y.
\end{aligned}
\end{equation}
Therefore,   the eigenvectors $v=(v_1,  v_2,  .  .  .  ,  v_{N_1})^{T}$ associated with $\lambda=\frac{2}{n-1}$ can be determined by Eqs.  (\ref{vk'.})(\ref{vivj})  and (\ref{xevit}).   Notice that
$v_{k^e_1}=v_{k^e_2}=\cdots=v_{k^e_{n-2}}=x_e$.  Let $\mbox{\bf x}=(x_e)^T$ which is an $E_0$ dimensional vector by the construction of $CL(G)$.

For convenience, we assume that the first $|X|$ rows of the incident matrix $B$ of $G$ correspond to the vertices of $X$, and hence the matrix $B$ can be written as $B={B_X \choose B_Y}$.   Let $D_X$ and $D_Y$ denote the
volume vectors which consist of degree sequences of vertices of $X$ and $Y$, respectively. Let
\begin{equation*}\label{bxdit}
C=\left(
\begin{array}{cc}
B_X&-D_X \\
B_Y&D_Y
\end{array}
\right).
\end{equation*}
Hence Eqs.  (\ref{vk'.})(\ref{vivj})  and (\ref{xevit}) are equivalent to the equation system $C {\mbox{\bf x} \choose t}=0$.

By Lemma \ref{bx}, the rank of $B$ is $N_0-1$ when $G$ is bipartite. Now we need to determine the rank of $C$. We denote the volume vectors of $C$ by $\mbox{\bf e}_1, \mbox{\bf e}_2, \cdots, \mbox{\bf e}_{E_0},\mbox{\bf
e}_0$ from left to right.
Assume that $\mbox{\bf e}_0$ is linearly related to the   $\mbox{\bf e}_1, \mbox{\bf e}_2, \cdots, \mbox{\bf e}_{E_0}$, it means that, there exist  constants $c_1, c_2, \cdots, c_{E_0}$  making the following formula true,
\begin{equation}\label{et}
 \mbox{\bf e}_0= c_1 \mbox{\bf e}_1+ c_2 \mbox{\bf e}_2+ \cdots+ c_{E_0} \mbox{\bf e}_{E_0}.
\end{equation}
 For every volume of $C$, there are two entries $1$ in $B_X$ and $B_Y$, respectively. From Eq. \eqref{et}, by summing all the first $|X|$ entries in $\mbox{\bf e}_0$, we have $c_1 + c_2 + \cdots + c_{E_0}= \sum\limits_{i=1}^{|X|} (-d_i)$.  For the same reason, we can get $c_1 + c_2 + \cdots + c_{E_0}=
 \sum\limits_{i=|X|+1}^{N_0} d_i$.  This implies that $\sum\limits_{i=1}^{|X|} (-d_i)=\sum\limits_{i=|X|+1}^{N_0} d_i$. Notice that $d_i>0$ for each $i=1,2,\cdots,N_0$. Hence it is obvious that $\sum\limits_{i=1}^{|X|}
 (-d_i)= \sum\limits_{i=|X|+1}^{N_0} d_i$ is impossible.  Thus we get a contradiction.   So $\mbox{\bf e}_0$ and $\mbox{\bf e}_1, \mbox{\bf e}_2, \cdots, \mbox{\bf e}_{E_0}$  are  linearly independent, i.e., the rank of matrix $C$ is
 $r(C)=r(B)+1=N_0$.

 Therefore, the number of solutions in  basic solution system of $C {\mbox{\bf x} \choose t}=0$ is $E_0-N_0+1$ when $G$ is bipartite,   i.e.,   $m_{\mathcal{L}_C}(\frac{2}{n-1})=E_0-N_0+1$.\vskip2mm

{\upshape (\romannumeral5)} Substituting $\lambda=\frac{n}{n-1}$ into Eq. \eqref{vk1}, we have
\begin{equation*}
v_{k^e_1}+v_{k^e_2}+v_{k^e_3}+\cdots+v_{k^e_{n-2}}+\frac{v_{i}}{\sqrt{d_{i}}}+\frac{v_{j}}{\sqrt{d_{j}}}=0.
\end{equation*}

For convenience, for each edge $e_s\in E(G)$, $s=1, 2, \ldots, E_0$, denote by $i_s$ and $j_s$ the end vertices of $e_s$.
So, we have the following linear equation system
\begin{equation}\label{last}
\begin{cases}
\begin{aligned}
&v_{k^{e_1}_1}&+&v_{k^{e_1}_2}&+&v_{k^{e_1}_3}&+&\cdots&+&v_{k^{e_1}_{n-2}}&+&\frac{v^{e_1}_{i_1}}{\sqrt{d^{e_1}_{i_1}}}&+&\frac{v^{e_1}_{j_1}}{\sqrt{d^{e_1}_{j_1}}}&=0,\\
&v_{k^{e_2}_1} &+&v_{k^{e_2}_2}&+&v_{k^{e_2}_3}&+&\cdots&+&v_{k^{e_2}_{n-2}}&+&\frac{v^{e_2}_{i_2}}{\sqrt{d^{e_2}_{i_2}}} &+&\frac{v^{e_2}_{j_2}}
  {\sqrt{d^{e_2}_{j_2}}}&=0,\\
&\vdots & & & & & & & & & & & & \\
&v_{k^{e_{E_0}}_1}&+& v_{k^{e_{E_0}}_2}&+&v_{k^{e_{E_0}}_3}&+&\cdots&+&v_{k^{e_{E_0}}_{n-2}}&+&\frac{v^{e_{E_0}}_{i_{E_0}}}
  {\sqrt{d^{e_{E_0}}_{i_{E_0}}}}&+&\frac{v^{e_{E_0}}_{j_{E_0}}}
  {\sqrt{d^{e_{E_0}}_{j_{E_0}}}}&= 0.
\end{aligned}
\end{cases}
\end{equation}
The corresponding coefficient matrix contains the following $E_0\times (n-2)E_0$ submatrix
\begin{equation*}
\left(
\begin{array}{cc}
{\underbrace{1\quad \cdots \quad 1}_{n-2} \quad { 0\quad \cdots \quad 0}\quad \cdots \quad {0\quad \cdots \quad 0}} \\
{{0\quad \cdots \quad 0} \quad\underbrace{ 1\quad \cdots \quad 1}_{n-2}\quad \cdots \quad { 0\quad \cdots \quad 0}} \\
\vdots \\
{{0\quad \cdots \quad 0} \quad { 0\quad \cdots \quad 0}\quad \cdots \quad\underbrace{ 1\quad \cdots \quad 1}_{n-2}}
\end{array}
\right).
\end{equation*}
Clearly, the submatrix above is of rank $E_0$.  Hence the number of solutions in a basic solution system of the system (\ref{last}) is $(n-3)E_0+N_0$. Therefore, $m_{\mathcal{L}_C}(\frac{n}{n-1})=(n-3)E_0+N_0$.\\
This completes the  proof of the theorem.
\end{pf}

\section{Related indexes and clique-blew up iterative graph}\label{applications}

Let $CL_0(G)=G$ and $CL_r(G)=CL(CL_{r-1}(G))$ for $r\geq 1$. The graph $CL_r(G)$ is called the $r$-th clique-blew up iterative graph of $G$. The number of vertices and edges of $CL_r(G)$, $r\geq0$, are denoted by $N_r$ and $E_r$, respectively. From the iterative method of  the  clique-blew up graph, we have
$$E_r=\frac{n(n-1)E_{r-1}}{2} \ \mbox{and}\  N_r=N_{r-1}+(n-2)E_{r-1}.$$
Hence
\begin{equation}\label{E_rN_r}
\quad E_r=\frac{n^{r}\cdot (n-1)^{r}E_{0}}{2^r}  \ \mbox{and}\ N_r=N_{0}+\frac{2E_{0}(\frac{n^r\cdot (n-1)^r}{2^r}-1)}{n+1}.
\end{equation}
For convenience, denote  by $\mathcal L_r$ the normalized Laplacian of $CL_r(G)$ for $r\geq 0$. Denote by $\sigma_r$ the normalized Laplacian spectrum of $CL_r(G)$ for $r\geq 0$.   From Theorem \ref{sp-Q}, we have the
following theorem.
\begin{thm}\label{sigm}
Let $G$ be a simple connected graph. For  $r\geq 2$ and $n\geq 3$,
$$\sigma_r=
\{{x\over n-1}| x\in \sigma_{r-1}\backslash \{0\}\}\cup\{0,\frac{2}{n-1},\frac{n}{n-1}\},$$
where
$m_{{\mathcal{L}_r}}(\frac{x}{n-1})=m_{{\mathcal{L}_{r-1}}}(x)$ for $x\in \sigma_{r-1}\backslash \{0\}$, $m_{{\mathcal{L}_r}}(0)=1$, $m_{{\mathcal{L}_r}}(\frac{2}{n-1})=E_{r-1}-N_{r-1}$ and  $m_{{\mathcal{L}_r}}(\frac{n}{n-1})=(n-3)E_{r-1}+N_{r-1}$.
\end{thm}

\begin{thm}
Let $G$ be a simple connected graph. For  $r\geq 1$ and $n\geq 3$, the multiplicative degree-Kirchhoff index $Kf^{*}(CL_r(G)$ of the $r$-clique-blew up graph $CL_r(G)$ can be determined by the multiplicative degree-Kirchhoff index $Kf^{*}(G)$ of the initial graph $G$ as follows
\begin{small}
\begin{equation}\label{KfQ_r}
\begin{aligned}
Kf^{*}(CL_r(G))=&\frac{n^r\cdot (n-1)^{2r}}{2^r}Kf^{*}(G)-
\frac{n^{r-1}\cdot (n-1)^{2r+1}}{2^r}(1-\frac{1}{(n-1)^r})E_0N_0\\
&+\frac{(n-1)^{2r}\cdot n^{r-1}}{2^{r-1}}\left(3(\frac{n^r}{2^r}-1)-\frac{1}{n+1}
(\frac{n^r}{2^{r-1}}+\frac{1}{(n-1)^{r-1}}-n-1)\right){E_0}^2.
\end{aligned}
\end{equation}
\end{small}
\end{thm}
\begin{pf}
Recall the normalized Laplacian eigenvalues of $G$ is $0=\lambda_1<\lambda_2\leq\cdots\leq\lambda_{N_0}$. Whether $G$ is bipartite or not, we have the following result by Theorem \ref{sp-Q} and Lemma \ref{lm1.2}
\upshape (\romannumeral1)
\begin{small}
\begin{equation}\label{KfQ}
\begin{aligned}
Kf^{*}(CL(G))&=2E_1\left(\sum\limits_{i=2}^{N_0}\frac{n-1}{\lambda_i}+
\frac{n-1}{2}(E_0-N_0)+\frac{n-1}{n}((n-3)E_0+N_0)\right)\\
&=\frac{1}{2}n(n-1)^2Kf^{*}(G)+\frac{3}{2}(n-1)^2(n-2){E_0}^2-
\frac{1}{2}(n-1)^2(n-2)E_0N_0.
\end{aligned}
\end{equation}
\end{small}
From Eqs. (\ref{E_rN_r}) and (\ref{KfQ}) and the definition of the $r$-th clique-blew up iterative graph, we can get
\begin{small}
\begin{equation*}\label{KfQr,r-1}
\begin{aligned}
Kf^{*}(CL_r(G))&=\frac{1}{2}n(n-1)^2Kf^{*}(CL_{r-1}(G))+
\frac{3}{2}(n-1)^2(n-2){E_{r-1}}^2\\
&\quad -\frac{1}{2}(n-1)^2(n-2)E_{r-1}N_{r-1}\\
&=\frac{n^r\cdot (n-1)^{2r}}{2^r}Kf^{*}(G)-
\frac{n^{r-1}\cdot (n-1)^{2r+1}}{2^r}(1-\frac{1}{(n-1)^r})E_0N_0\\
&\quad+\frac{(n-1)^{2r}\cdot n^{r-1}}{2^{r-1}}\left(3(\frac{n^r}{2^r}-1)-\frac{1}{n+1}
(\frac{n^r}{2^{r-1}}+\frac{1}{(n-1)^{r-1}}-n-1)\right){E_0}^2.
\end{aligned}
\end{equation*}
\end{small}
The proof is completed.
\end{pf}
\begin{thm}
For  $r\geq 1$ and $n\geq 3$, the Kemeny's constant $K_e(CL_r(G))$ for the  random walks on $CL_r(G)$ is as follows
\begin{small}
\begin{equation*}
\begin{aligned}
K_e(CL_r(G))=&(n-1)^rK_e(G)+\frac{(n-1)^{r+1}}{2n}(\frac{1}{(n-1)^r}-1)N_0+\\
&\left(\frac{3(n-1)^r}{n}(\frac{n^r}{2^r}-1)+\frac{(n-1)^r}
{n+1}(1-\frac{n^{r-1}}{2^{r-1}})+\frac{(n-1)^{r-1}}
{n(n+1)}(1-\frac{1}{(n-1)^{r-1}})\right)E_0.
\end{aligned}
\end{equation*}
\end{small}
\end{thm}
\begin{pf}
By Lemma \ref{lm1.2} \upshape (\romannumeral4) and Eq. (\ref{KfQ}), it follows that
\begin{equation}\label{kecl}
\begin{aligned}
K_e(CL(G))&=\frac{1}{2E_1}Kf^{*}(CL(G))\\
&=\frac{n-1}{2E_0}Kf^{*}(G)+\frac{3(n-1)(n-2)}{2n}E_0-
\frac{(n-1)(n-2)}{2n}N_0\\
&=(n-1)K_e(G)+\frac{3(n-1)(n-2)}{2n}E_0-\frac{(n-1)(n-2)}{2n}N_0.
\end{aligned}
\end{equation}
From Eqs. (\ref{E_rN_r}) and (\ref{kecl}) and the definition of the $r$-th clique-blew up iterative graph, we can get
\begin{small}
\begin{equation*}
\begin{aligned}
K_e(CL_r(G))&=(n-1)K_e(CL_{r-1}(G))+\frac{3(n-1)(n-2)}
{2n}E_{r-1}-\frac{(n-1)(n-2)}{2n}N_{r-1}\\
&=(n-1)^rK_e(G)+\frac{(n-1)^{r+1}}{2n}(\frac{1}{(n-1)^r}-1)N_0+\\
&\quad\left(\frac{3(n-1)^r}{n}(\frac{n^r}{2^r}-1)+\frac{(n-1)^r}
{n+1}(1-\frac{n^{r-1}}{2^{r-1}})+\frac{(n-1)^{r-1}}
{n(n+1)}(1-\frac{1}{(n-1)^{r-1}})\right)E_0.
\end{aligned}
\end{equation*}
\end{small}
The proof is completed.
\end{pf}
\begin{thm}
For  $r\geq 1$ and $n\geq 3$, the number of spanning trees of $CL_r(G)$ is as follows
\begin{equation*}
\tau(CL_r(G))=2^{2E_0\alpha-rN_0-\frac{2E_0}{n+1}( 2\alpha-r )+r}\cdot n^{2(n-3)E_0\alpha+rN_0+\frac{2E_0}{n+1}(2\alpha-r)-r}\cdot \tau(G),
\end{equation*}
where  $\alpha=\frac{\frac{n^r\cdot (n-1)^r}{2^r}-1}{n^2-n-2} $.

\end{thm}
\begin{pf}
Let the normalized Laplacian eigenvalues of $CL(G)$ be $0=\lambda^{'}_1 < \lambda^{'}_2 \leq\cdots\leq\lambda^{'}_{N_1}$. Whether $G$ is bipartite or not, by Lemma \ref{lm1.2} \upshape
(\romannumeral3) and the definition of $CL(G)$ we have
\begin{equation}\label{tauQ/tauG}
\frac{\tau(CL(G))}{\tau(G)}=\frac{2(n-1)^{N_0+
(n-2)E_0-1}\cdot \prod\limits_{i=2}^{N_1}\lambda^{'}_{i}}
{n\prod\limits_{i=2}^{N_0}\lambda_i} .
\end{equation}
From Theorem \ref{sp-Q} we have
\begin{equation}\label{mmlambda}
\begin{aligned}
\prod\limits_{i=2}^{N_1}\lambda^{'}_{i}
&=(\frac{1}{n-1})^{N_0-1}\cdot (\frac{2}{n-1})^{E_0-N_0}\cdot
 (\frac{n}{n-1})^{(n-3)E_0+N_0}\cdot \prod\limits_{i=2}^{N_0}\lambda_i\\
&=\frac{2^{E_0-N_0}\cdot n^{(n-3)E_0+N_0}}{(n-1)^{N_0+
  (n-2)E_0-1}}\prod\limits_{i=2}^{N_0}\lambda_i.
\end{aligned}
\end{equation}
By Eqs. (\ref{tauQ/tauG}) and (\ref{mmlambda}), we have
\begin{equation*}\label{tauQ}
\tau(CL(G))=2^{E_0-N_0+1}\cdot n^{(n-3)E_0+N_0-1}\cdot \tau(G).
\end{equation*}
It follows from the recursive relation that
\begin{small}
\begin{equation*}
\begin{aligned}
\tau(CL_r(G))&=2^{E_{r-1}-N_{r-1}+1}\cdot n^{(n-3)E_{r-1}+N_{r-1}-1}\cdot \tau(CL_{r-1}(G))\\
&=2^{{\sum\limits_{i=0}^{r-1}(E_i-N_i)}+r}\cdot
  n^{{\sum\limits_{i=0}^{r-1}((n-3)E_i-N_i)}-r}\cdot \tau(G)\\
&=2^{2E_0\alpha-rN_0-\frac{2E_0}{n+1}( 2\alpha-r )+r}\cdot n^{2(n-3)E_0\alpha+rN_0+\frac{2E_0}{n+1}(2\alpha-r)-r}\cdot \tau(G).
\end{aligned}
\end{equation*}
\end{small}
The proof is completed.
\end{pf}

\vskip 1cm \noindent {\bf Acknowledgement} This work was supported by National Natural Science Foundation of China (11571320 and 11671366) and  Zhejiang Provincial Natural Science Foundation (LY19A010018).



\end{CJK}
\end{document}